\documentclass{elsart3-1}

\usepackage{amssymb}

\usepackage[english,francais]{babel}

\newtheorem{theorem}{Theorem}[section]
\newtheorem{lemma}[theorem]{Lemma}
\newtheorem{e-proposition}[theorem]{Proposition}

\newtheorem{e-definition}[theorem]{Definition\rm}
\newtheorem{remark}{\it Remark\/}

\renewcommand{\theequation}{\arabic{equation}}
\def\og{\leavevmode\raise.3ex\hbox{$\scriptscriptstyle\langle\!\langle$~}}
\def\fg{\leavevmode\raise.3ex\hbox{~$\!\scriptscriptstyle\,\rangle\!\rangle$}}

\journal{the Acad\'emie des sciences}
\begin{document}
\centerline{}
\begin{frontmatter}


\selectlanguage{english}
\title{A generalization of Cram\'{e}r large deviations for martingales}


\selectlanguage{english}
\author[label1,label2]{Xiequan Fan}
\ead{fanxiequan@hotmail.com}
\author[label1]{Ion Grama}
\ead{ion.grama@univ-ubs.fr}
\author[label1]{Quansheng Liu}
\ead{quansheng.liu@univ-ubs.fr}
\address[label1]{Univ. Bretagne-Sud,  UMR 6205, LMBA, F-56000 Vannes, France}
\address[label2]{Regularity Team, Inria and MAS Laboratory, Ecole Centrale Paris - Grande Voie des Vignes,\\ 92295 Ch\^{a}tenay-Malabry, France}


\medskip
\begin{center}
{\small Received *****; accepted after revision +++++\\
Presented by £££££}
\end{center}

\begin{abstract}
\selectlanguage{english}
In this note, we give a generalization of Cram\'{e}r's large deviations for
martingales, which  can be regarded as a supplement of  Fan,  Grama and Liu  (Stochastic Process. Appl., 2013).
Our method is based on the change of probability measure developed by
Grama and Haeusler (Stochastic Process. Appl., 2000).
{\it To cite this article: X. Fan, I. Grama,  Q. Liu,  A generalization of Cram\'{e}r large deviations for martingales, C. R. Acad. Sci. Paris, Ser. I 352 (2014) 853--858. }

\vskip 0.5\baselineskip

\selectlanguage{francais}
\noindent{\bf R\'esum\'e} \vskip 0.5\baselineskip \noindent
{\bf Une g\'{e}n\'{e}ralisation des grandes d\'{e}viations de Cram\'{e}r  pour les martingales.}
Dans cette note, nous donnons une g\'{e}n\'{e}ralisation des grandes d\'{e}viations de Cram\'{e}r pour les martingales, qui peut \^{e}tre consid\'{e}r\'{e} comme un  suppl\'{e}ment de Fan, Grama et Liu (Stochastic Process. Appl., 2013). Notre m\'{e}thode est bas\'{e}e sur le changement de mesure de probabilit\'{e} d\'{e}velopp\'{e} par Grama et Haeusler (Stochastic Process. Appl., 2000). 
\end{abstract}
\end{frontmatter}

\section{Introduction}
Assume that $\eta_{1},..., \eta_{n}$ is a sequence of independent and identically distributed (i.i.d.) centered real valued random variables  satisfying the following Cram\'{e}r condition: $\mathbb{E}\exp\{ c_{0}|\eta_{1}|\}<\infty$ for some $c_{0}>0$. Denote by $\sigma^2=\mathbb{E}\eta_{1}^2, \xi_{i}=\eta_i/(\sqrt{n} \sigma)$ and $X_n=\sum_{i=1}^{n}\xi_{i}.$  Cram\'{e}r \cite{Cramer38} has established the following asymptotic expansion of the tail probabilities of $X_n$, for all $0 \leq x =o( n^{1/6})$ as $n\rightarrow \infty,$
\begin{equation}
\mathbb{P}(X_n> x )=\Big(1-\Phi(x)\Big)\Big(1+o(1)\Big),
\label{cramer001}
\end{equation}
where
\[
\Phi(x)=\frac{1}{\sqrt{2\pi}}  \int_{-\infty}^{x}\exp\Big\{-\frac{t^2}{2} \Big\} dt
\]
 is the standard normal distribution function.
More precise results can be found in Feller \cite{Fl43}, Petrov \cite{Pe54,Petrov75},  Saulis and Statulevi\v{c}ius \cite{SS78},   Sakhanenko \cite{S91} and \cite{F12} among others.

Let $(\xi _i,\mathcal{F}_i)_{i=0,...,n}$ be a sequence of martingale differences  defined on some
 probability space $(\Omega ,\mathcal{F},\mathbb{P})$,  where $\xi
_0=0 $ and $\{\emptyset, \Omega\}=\mathcal{F}_0\subseteq ...\subseteq \mathcal{F}_n\subseteq
\mathcal{F}$ are increasing $\sigma$-fields. Set
\begin{equation}
X_{0}=0,\ \ \ \ \ X_k=\sum_{i=1}^k\xi _i,\quad k=1,...,n.  \label{xk}
\end{equation}
Denote by $\left\langle X\right\rangle $ the quadratic characteristic of the
martingale $X=(X_k,\mathcal{F}_k)_{k=0,...,n}:$
\begin{equation}\label{quad}
\left\langle X\right\rangle _0=0,\ \ \ \ \ \left\langle X\right\rangle _k=\sum_{i=1}^k\mathbb{E}(\xi _i^2|\mathcal{F}
_{i-1}),\quad k=1,...,n.
\end{equation}
Consider the stationary case for simplicity. For the martingale differences having a $(2+p)$th moment, i.e.  $||\xi_i||_{2+p}< \infty$ for some $p \in (0, 1],$  expansions of the type (\ref{cramer001}) in the range $0 \leq x =o(\sqrt{\log n} )$ have been obtained by Haeusler and Joos \cite{HJ88}, Grama \cite{G97} and  Grama and Haeusler \cite{GH06}.  If the martingale differences  are bounded $|\xi_i |\leq C/\sqrt{n} $  and satisfy $ ||\left\langle X\right\rangle _n-1  ||_{\infty} \leq L^2/n \ a.s.\ $for two positive  constants $C$ and  $L$,
expansion (\ref{cramer001})  has been firstly established by Ra\v{c}kauskas \cite{Rackauskas90,Rackauskas95} in the range $0 \leq x =o(n^{1/6}  )$, and then this range has been extended to a larger one $0 \leq x =o( n^{1/4}  )$ by Grama and Haeusler \cite{GH00} with a method based on change of probability measure.
Recently, Fan et al.\ \cite{F13} have generalized the result of Grama and Haeusler \cite{GH00} to a much larger range $0 \leq x =o( n^{1/2})$ for  $\xi_i$  satisfying the following conditional Bernstein  condition: for a positive   constant $C,$
\begin{equation}
|\mathbb{E}(\xi_{i}^{k}  | \mathcal{F}_{i-1})| \leq \frac{1}{2}\, k!\,\Big(\frac{C}{\sqrt{n}} \Big)^{k-2} \mathbb{E}(\xi_{i}^2 | \mathcal{F}_{i-1})
\ \ \  \mbox{for all} \  k\geq 2 \  \textrm{and all}   \ 1\leq i\leq n.
\label{Bernst cond}
\end{equation}
It is worth noting that the conditional Bernstein  condition does not imply that $\xi_i$'s are bounded.

The aim of this note is to extend the expansion of Fan et al.\  \cite{F13}  to the case of martingale differences satisfying  the following conditional Cram\'{e}r condition considered in Liu and Watbled \cite{Liu09a}:
\begin{equation}
\;\sup_{ i } \mathbb{E} (\exp\{C_0 \sqrt{n} |\xi_{i}|\} |\mathcal{F}_{i-1}) \leq C_{1}\, ,\label{cra cond}
\end{equation}
where $C_0$ and $C_1$ are two positive constants. It is worth noting that,  in general, condition (\ref{cra cond}) does not imply the conditional Bernstein  condition (\ref{Bernst cond}), unless
$n\mathbb{E}(\xi_{i}^2 | \mathcal{F}_{i-1})$ are all bounded from below by a positive constant. Thus our result is not a  consequence of  Fan et al.\ \cite{F13}.

Throughout this paper, $c$ and $c_\alpha,$ probably supplied with some indices,
denote respectively a generic positive absolute constant and a generic positive constant depending only on $\alpha.$
Moreover, $\theta$ stands for any value satisfying $\left| \theta  \right| \leq 1$.

\section{\textbf{Main Results} \label{sec2}}

The following theorem is our main result, which  can be regarded as a parallel result of  Fan et al.\ \cite{F13}
under the conditional Cram\'{e}r condition:
\begin{description}
\item[(A1)] $\;\sup_{1\leq i\leq n} \mathbb{E} (\exp\{c_0n^{1/2}|\xi_{i}|\} |\mathcal{F}_{i-1}) \leq c_{1}\, $;
\item[(A2)]\ \  $ ||\left\langle X\right\rangle _n-1 ||_{\infty} \leq \delta^2 \ \  a.s.,\   \  $ where $\delta$ is nonnegative and usually depends on $n.$
\end{description}
\begin{theorem}\label{th1}
Assume conditions (A1) and (A2). Then there
exists a positive absolute constant $\alpha_0,$ such that
for all $0\leq x \leq \alpha_0 \, n^{1/2}$ and $\delta \leq \alpha_0$, the following equalities hold
\begin{equation} \label{t1ie1}
\frac{\mathbb{P}(X_n>x)}{1-\Phi \left( x\right)} = \exp\Bigg\{\theta c_{\alpha_0} \! \Bigg( \frac{ x^3}{\sqrt{n}} + x^2 \delta^2 + (1+ x) \Big( \frac{  \log n}{\sqrt{n}} +  \delta \Big)\Bigg) \Bigg\}
\end{equation}
and
\begin{equation}
  \frac{\mathbb{P}(X_n<-x)}{\Phi \left( -x\right)} = \exp\Bigg\{\theta c_{\alpha_0} \! \Bigg( \frac{ x^3}{\sqrt{n}} + x^2 \delta^2 + (1+ x) \Big( \frac{  \log n}{\sqrt{n}} +  \delta \Big)\Bigg) \Bigg\},
\end{equation}
where $|\theta|\leq 1$.
In particular, for all $0\leq x =o \big( \min\{n^{1/6},
\delta^{-1}\}\big)$ as $\min\{n, \delta^{-1}\} \rightarrow \infty$,
\begin{equation}
 \mathbb{P}(X_n\geq x) =\Big(1- \Phi \left( x\right)\Big)\Big(1 + o(1)\Big).
\end{equation}
\end{theorem}

From  (\ref{t1ie1}), we find that there is an absolute constant $\alpha_0 >0$ such that
for all $0\leq x \leq \alpha_0  \, n^{1/2}$ and $\delta \leq \alpha_0$,
\begin{equation}\label{fsvdf1}
\Bigg| \log \frac{\mathbb{P}(X_n>x)}{1-\Phi \left( x\right)} \Bigg|\  \leq  \     c_{\alpha_0} \! \Bigg( \frac{ x^3}{\sqrt{n}} + x^2 \delta^2 + (1+ x) \Big( \frac{  \log n}{\sqrt{n}} +  \delta \Big)\Bigg).
\end{equation}
Note that this result can be regarded as a refinement of the moderate deviation principle (MDP) under  conditions (A1) and (A2).
Let $a_n$ be any sequence of real numbers satisfying $a_n \rightarrow \infty$ and $a_nn^{-1/2}\rightarrow 0$
as $n\rightarrow \infty$. If $\delta\rightarrow 0$ as $n\rightarrow \infty,$
then inequality (\ref{fsvdf1}) implies the MDP for $X_n$ with the speed $a_n$ and good rate function $x^2/2;$ for each Borel set $B$,
\begin{eqnarray*}
- \inf_{x \in B^o}\frac{x^2}{2} \ \leq \ \liminf_{n\rightarrow \infty}\frac{1}{a_n^2}\log \mathbb{P}\left(\frac{1}{a_n} X_n \in B \right)
\ \leq \ \limsup_{n\rightarrow \infty}\frac{1}{a_n^2}\log \mathbb{P}\left(\frac{1}{a_n} X_n \in B \right) \leq  - \inf_{x \in \overline{B}}\frac{x^2}{2} \, ,
\end{eqnarray*}
where $B^o$ and $\overline{B}$ denote the interior and the closure of $B$, respectively (see Fan et al.\ \cite{F13} for   details).

\section{\textbf{Sketch of the proof}\label{sec3}}
Let $(\xi _i,\mathcal{F}_i)_{i=0,...,n}$ be a martingale differences satisfying the condition (A1).
 For any real number $\lambda$ with $|\lambda| \leq c_0 n^{1/2} ,$ define
\[
Z_k(\lambda )=\prod_{i=1}^k\frac{e^{\lambda \xi _i}}{\mathbb{E}(e^{\lambda \xi _i}|
\mathcal{F}_{i-1})},\quad k=1,...,n,\quad Z_0(\lambda )=1.  \label{C-1}
\]
Then $Z(\lambda )=(Z_k(\lambda ),\mathcal{F}_k)_{k=0,...,n}$  is a positive martingale and  for each
real number $\lambda$ with $|\lambda| \leq c_0 n^{1/2}$ and each $k=1,...,n,$ the random variable $Z_k(\lambda
) $ is a probability density on $(\Omega ,\mathcal{F},\mathbb{P}).$ Thus we can
 define the \emph{conjugate probability measure} $\mathbb{P}_\lambda $ on $(\Omega ,%
\mathcal{F})$, where
\begin{equation}
d\mathbb{P}_\lambda =Z_n(\lambda )d\mathbb{P}.  \label{f21}
\end{equation}
Denote by $\mathbb{E}_{\lambda}$ the expectation with respect to $\mathbb{P}_{\lambda}.$ Setting
\[
b_i(\lambda)=\mathbb{E}_{\lambda}(\xi_i |\mathcal{F}_{i-1})   \quad  \textrm{and}\quad  \eta_i(\lambda)=\xi_i - b_i(\lambda)\ \  \ \textrm{for}\  i=1,...,n,
\]
we obtain the decomposition of $X_n$ similar to that of Grama and Haeusler \cite{GH00}:
\begin{equation}
X_n=B_n(\lambda )+Y_n(\lambda ), \label{xb}
\end{equation}
where
\[
B_n(\lambda )=\sum_{i=1}^n b_i(\lambda )\ \ \ \ \textrm{and}\ \ \ \ \ Y_n(\lambda )=\sum_{i=1}^n\eta _i(\lambda ).
\]
Note that $(Y_k( \lambda  ), \mathcal{F}_{k})_{k=1,...,n}$ is also a  sequence of martingale differences w.r.t. $\mathbb{P}_{\lambda}$.
\vspace{0.2cm}

In the sequel, we establish some auxiliary lemmas which will be used in the proof of Theorem \ref{th1}.
We first give upper bounds for the conditional moments.
\begin{lemma}\label{lema1}  Assume condition (A1). Then
\[
\mathbb{E}(|\xi_i|^k | \mathcal{F}_{i-1}) \leq    k!\,(c_0n^{1/2})^{-k}c_1, \ \ \ \  \ k\geq 3.
\]
\end{lemma}
\noindent  \emph{Proof. }
 Applying the elementary inequality $ x^k/ k! \leq  e^{x}$
to $x=c_0|n^{1/2}\xi_i|$, we have, for $k\geq 3,$
\begin{equation}\label{fghk2}
 |\xi_i|^k  \ \leq \ k!\,(c_0n^{1/2})^{-k} \exp\{c_0|n^{1/2}\xi_i|\} .
\end{equation}
Taking conditional expectations on both sides of the last inequality, by condition (A1), we obtain   the desired inequality.   \qed

\begin{remark}
It is worth noting that both  condition (A1) and  the conditional  Bernstein  condition (\ref{Bernst cond}) imply the following  hypothesis.
\begin{description}
\item[\emph{(A1$'$)}] There exists  $\epsilon > 0,$ usually depends on $n$, such that
\[
\mathbb{E}(|\xi_i|^k | \mathcal{F}_{i-1}) \leq   c_1 \, k!\, \epsilon^{k}   \ \  \ \mbox{for all}\  k\geq 2 \  \textrm{and all}   \ 1\leq i\leq n.
\]
\end{description}
When $\epsilon=c_2/ \sqrt{n},$ condition (A1$'$), together (A2), yields Theorem \ref{th1}.
\end{remark}

Using Lemma \ref{lema1}, we obtain the following two technical lemmas. Their proofs are similar to the arguments of Lemmas 4.2 and 4.3 of Fan et al.\ \cite{F13}.
\begin{lemma}
\label{lemma5} Assume conditions (A1) and (A2). Then, for all $0 \leq \lambda \leq \frac14 c_0 n^{1/2} ,$
\begin{equation}\label{f82}
\left| B_n(\lambda )-\lambda \right| \ \leq \ c\, (\lambda \delta^2 + \lambda^2 n^{-1/2}).
\end{equation}
\end{lemma}

\begin{lemma}
\label{lemma6} Assume conditions (A1) and (A2). Then, for all $0 \leq \lambda \leq \frac14 \, c_0 n^{1/2} ,$
\[
\left| \Psi _n(\lambda )-\frac{\lambda ^2}2\right| \ \leq \  c\, ( \lambda^2 \delta^2  +\lambda ^3n^{-1/2}) ,
\]
where \[\Psi _n(\lambda )=\sum_{i=1}^n \log \mathbb{E}(e^{\lambda \xi _i}|\mathcal{F}_{i-1}).\]
\end{lemma}

The following lemma gives the rate of convergence in the central limit theorem for the conjugate martingale $(Y_i( \lambda  ), \mathcal{F}_{i})$ under the probability measure $\mathbb{P}_{ \lambda  }.$ Its proof is similar to that of Lemma 3.1 of Fan et al.\ \cite{F13} by noting the fact that $\mathbb{E}(\xi_i^2 | \mathcal{F}_{i-1}) \leq c/n.$

\begin{lemma}
\label{LEMMA4}
Assume conditions (A1) and (A2). Then,  for all $0 \leq \lambda \leq \frac14 \, c_0 n^{1/2} ,$
\[
\sup_{x}\Big| \mathbb{P}_\lambda (\, Y_n(\lambda )\leq x)-\Phi (x)\Big| \ \leq\
c\left( \lambda \,\frac{1}{ \sqrt{n}} +\frac{\log  n}{ \sqrt{n}}  +\delta \right) .
\]
\end{lemma}

\vspace{0.3cm}

\noindent\emph{Proof of Theorem \ref{th1}. }   The proof of Theorem 2.1  is similar to the arguments of Theorems 2.1 and 2.2 in Fan et al.\ \cite{F13} with
$\epsilon=\frac{c_0}{4 \sqrt{n}}$. However, instead of using Lemmas 4.2, 4.3 and 3.1 of \cite{F13}, we shall make use of Lemmas \ref{lemma5},  \ref{lemma6} and \ref{LEMMA4} respectively. \qed

\section*{Acknowledgements}
We thank the reviewer for his/her thorough review and highly appreciate the comments and
suggestions, which significantly contributed to improving the quality of the publication.
The work has been partially supported by the National Natural Science Foundation of China (Grants No. 11171044 and No. 11101039), and by Hunan Provincial Natural Science Foundation of China (Grant No. 11JJ2001).

\selectlanguage{english}

\newpage

The proofs of  Lemmas
\ref{lemma5} and \ref{lemma6} are given below.

\noindent \emph{Proof of Lemma
\ref{lemma5}. } Recall that $0 \leq \lambda \leq \frac14 \, c_0 n^{1/2} .$ By the relation between $\mathbb{E}$ and $\mathbb{E}_{\lambda}$ on $\mathcal{F}_i,$
we have
\[
b_i(\lambda )=\frac{\mathbb{E}(\xi
_ie^{\lambda \xi _i}|\mathcal{F}_{i-1})}{\mathbb{E}(e^{\lambda \xi _i}|\mathcal{F}%
_{i-1})},\quad\ \ \ \ \ \   i=1,...,n.
\]
Jensen's inequality and $\mathbb{E}(\xi _i|\mathcal{F}_{i-1})=0$ imply that $\mathbb{E}(e^{\lambda \xi _i}|\mathcal{F}
_{i-1})\geq 1.$ Since
\[
\mathbb{E}(\xi_{i} e^{\lambda\xi_{i}} |\mathcal{F}_{i-1})=\mathbb{E}\left(\xi_{i}(e^{\lambda\xi_{i}}-1)|\mathcal{F}_{i-1} \right)\geq 0,\
\]
 by Taylor's expansion for $e^x$, we find that
\begin{eqnarray}
B_n(\lambda ) & \leq & \sum_{i=1}^{n}\mathbb{E}(\xi_{i} e^{\lambda \xi_{i}} | \mathcal{F}_{i-1})  =    \lambda\langle X\rangle_{n}+ \sum_{i=1}^{n}\sum_{k=2}^{+\infty}\mathbb{E} \left(\frac{\xi_{i}(\lambda\xi_{i})^{k}}{k !}  \Bigg |  \mathcal{F}_{i-1}   \right)   .\label{f26}
\end{eqnarray}
Using Lemma \ref{lema1}, we obtain
\begin{eqnarray}
  \sum_{i=1}^{n}\sum_{k=2}^{+\infty}\Bigg| \mathbb{E}\left(\frac{\xi_{i}(\lambda\xi_{i})^{k}}{k !}  \Bigg |  \mathcal{F}_{i-1}   \right)\Bigg|
  & \ \leq \ & \sum_{i=1}^{n}\sum_{k=2}^{+\infty}|\mathbb{E}\left( \xi_{i}^{k+1}| \mathcal{F}_{i-1} \right)| \frac{\lambda^k}{k !} \nonumber \\
  & \leq & \sum_{i=1}^{n}\sum_{k=2}^{+\infty} c_1(k+1) \lambda^k (c_0 n^{1/2})^{-k-1} \nonumber \\
& \leq &c_2  \, \lambda^2 n^{-1/2}  .\label{f52}
\end{eqnarray}
Condition (A2) together with (\ref{f26}) and (\ref{f52}) imply the upper bound of $B_n(\lambda )$:
\[
B_n(\lambda ) \leq
\lambda +  \lambda \delta^{2} + c_2  \, \lambda^2 n^{-1/2} .
\]
Using  Lemma \ref{lema1}, we have
\begin{eqnarray}
\mathbb{E}\left(e^{\lambda \xi_{i}} | \mathcal{F}_{i-1} \right) & \ \leq \ & 1 + \sum_{k=2}^{+\infty}\left|\mathbb{E} \left(\frac{(\lambda\xi_{i})^k}{k !} \Bigg| \mathcal{F}_{i-1} \right)\right|  \nonumber \\
& \leq & 1+   \sum_{k=2}^{+\infty}c_1 \lambda^k  \left(c_0 n^{1/2}\right)^{-k}  \nonumber \\
& \leq & 1+ c_{3}\, \lambda^2 n^{-1} .  \label{sbopdf}
\end{eqnarray}
This inequality together with (\ref{f52}) and condition (A2) imply the lower bound of $B_{n}(\lambda)$:
\begin{eqnarray*}
  B_n(\lambda ) &\ \geq \ & \left(\sum_{i=1}^{n}\mathbb{E}(\xi_{i} e^{\lambda \xi_{i}} | \mathcal{F}_{i-1})\right)\Bigg( 1+ c_{3}\, \lambda^2n^{-1}\Bigg)^{-1}\nonumber\\
  & \geq &  \left(\lambda\langle X\rangle_{n} -  \sum_{i=1}^{n}\sum_{k=2}^{+\infty}\left|\mathbb{E}\left(\frac{\xi_{i}(\lambda\xi_{i})^{k}}{k !} \Bigg | \mathcal{F}_{i-1} \right) \right| \right)\Bigg( 1+ c_{3}\, \lambda^2n^{-1} \Bigg)^{-1} \\
  & \geq &  \bigg(\lambda -\lambda \delta^2 -  c_2  \, \lambda^2 n^{-1/2}  \bigg)\bigg( 1+ c_{3}\, \lambda^2n^{-1} \bigg)^{-1} \\
  & \geq &  \lambda - \lambda \delta^{2} -   c_4  \, \lambda^2 n^{-1/2}.
\end{eqnarray*}
The proof of Lemma \ref{lemma5} is finished.\hfill\qed

\noindent \emph{Proof of Lemma
\ref{lemma6}. } Recall that $0 \leq \lambda \leq \frac14 \, c_0 n^{1/2}.$ Since $\mathbb{E}(\xi _i|\mathcal{F}_{i-1})=0$, it is easy to see that
\[
\Psi _n(\lambda )=\sum_{i=1}^n\left( \log \mathbb{E}(e^{\lambda \xi _i}|\mathcal{F}%
_{i-1})-\lambda \mathbb{E}(\xi _i|\mathcal{F}_{i-1})-\frac{\lambda ^2}2\mathbb{E}(\xi _i^2|%
\mathcal{F}_{i-1})\right) +\frac{\lambda ^2}2\left\langle X\right\rangle _n.
\]
Using a two-term  Taylor's expansion of $\log(1+x), x\geq0$,  we obtain
\begin{eqnarray*}
 \Psi _n(\lambda ) - \frac{\lambda ^2}2\left\langle X\right\rangle _n
&=& \sum_{i=1}^n\left( \mathbb{E}(e^{\lambda \xi _i}|\mathcal{F}%
_{i-1})-1-\lambda \mathbb{E}(\xi _i|\mathcal{F}_{i-1})-\frac{\lambda ^2}2\mathbb{E}(\xi _i^2|%
\mathcal{F}_{i-1}) \right) \\
&& - \frac1{2 \, \Big(1+|\theta|\big(\mathbb{E}(e^{\lambda \xi _i}|\mathcal{F}%
_{i-1})-1 \big) \Big)^2} \sum_{i=1}^n \left(\frac{}{} \mathbb{E}(e^{\lambda \xi _i}|\mathcal{F}%
_{i-1})-1 \right)^2.
\end{eqnarray*}
Since $\mathbb{E}(e^{\lambda \xi _i}|\mathcal{F}_{i-1})\geq1$, we find that
\begin{eqnarray*}
 \left|\Psi _n(\lambda ) - \frac{\lambda ^2}2\left\langle X\right\rangle _n \right|
&\leq& \sum_{i=1}^n\left| \mathbb{E}(e^{\lambda \xi _i}|\mathcal{F}%
_{i-1})-1-\lambda \mathbb{E}(\xi _i|\mathcal{F}_{i-1})-\frac{\lambda ^2}2\mathbb{E}(\xi _i^2|%
\mathcal{F}_{i-1}) \right|\\
&& + \frac12 \sum_{i=1}^n \left(\frac{}{} \mathbb{E}(e^{\lambda \xi _i}|\mathcal{F}%
_{i-1})-1 \right)^2\\
&\leq& \sum_{i=1}^{n}\sum_{k=3}^{+\infty}\frac{\lambda^{k}}{k !}|\mathbb{E}(\xi_{i}^{k} |\mathcal{F}_{i-1})| + \frac12\sum_{i=1}^{n} \left(\sum_{k=2}^{+\infty}\frac{\lambda^k}{k !}|\mathbb{E}(\xi_{i}^k|\mathcal{F}_{i-1})|\right)^2.
\end{eqnarray*}
In the same way as in the proof of Lemma \ref{lemma5}, by Lemma \ref{lema1},  we have
\[
\left|\Psi _n(\lambda ) - \frac{\lambda ^2}2\left\langle X\right\rangle _n \right|  \leq  c_3 \lambda^3 n^{-1/2} .
\]
Combining this inequality with condition (A2), we obtain the desired inequality.\hfill\qed

\begin{thebibliography}{00}{\footnotesize
\bibitem{Cramer38} Cram\'{e}r, H., 1938. Sur un nouveau th\'{e}or\`{e}me-limite de la
th\'{e}orie des probabilit\'{e}s. \textit{Actualite's Sci. Indust.}, \textbf{736}, 5--23.
\bibitem{F12} Fan, X., Grama, I. and Liu, Q., 2013. Sharp large deviations under Bernstein's condition. \emph{C. R. Acad. Sci. Paris, Ser. I} \textbf{351}, 845--848.
\bibitem{F13} Fan, X., Grama, I. and Liu, Q., 2013. Cram\'{e}r large deviation expansions for martingales under Bernstein's condition. \emph{Stochastic Process. Appl.} \textbf{123}, 3919--3942.
\bibitem{Fl43} Feller, W., 1943. Generalization of a probability limit theorem of Cram\'{e}r. \textit{Trans. Amer. Math. Soc.}, 361--372.
\bibitem{G97} Grama, I., 1997. On moderate deviations for
martingales. \emph{Ann. Probab.} \textbf{25}, 152--184.
\bibitem{GH00} Grama, I. and Haeusler, E., 2000. Large deviations
for martingales via Cram\'{e}r's method. \textit{Stochastic Process. Appl.} \textbf{85}, 279--293.
\bibitem{GH06}  Grama, I. and Haeusler, E.,  2006. An asymptotic expansion for
probabilities of moderate deviations for multivariate martingales. \textit{J.  Theoret. Probab.} \textbf{19}, 1--44.
\bibitem{HJ88} Haeusler, E. and Joos, K., 1988. A nonuniform bound on the rate of convergence in the martingale central limit theorem. \emph{Ann. Probab.} \textbf{16}, No. 4,    1699--1720.
\bibitem{Liu09a} Liu, Q. and Watbled, F., 2009. Exponential
ineqalities for martingales and asymptotic properties of the free energy of
directed polymers in a random environment. \emph{Stochastic Process. Appl.}
\textbf{119}, 3101--3132.

\bibitem{Pe54} Petrov, V. V., 1954. A generalization of Cram\'{e}r's limit theorem. \textit{Uspekhi Math. Nauk} \textbf{9}, 195--202.
\bibitem{Petrov75} Petrov, V. V., 1975. \textit{Sums of Independent
Random Variables.} Springer-Verlag. Berlin.
\bibitem{Rackauskas90}  Ra\v{c}kauskas, A., 1990. On probabilities of large
deviations for martingales. \textit{Liet. Mat. Rink.} \textbf{30},
784--795.
\bibitem{Rackauskas95}  Ra\v{c}kauskas, A., 1995. Large deviations for
martingales with some applications. \textit{Acta Appl. Math.}
\textbf{38}, 109--129.
\bibitem{S91}  Sakhanenko, A. I.,  1991. Berry-Esseen type bounds for large deviation probabilities, \textit{Siberian Math. J.} \textbf{32}, 647--656.
\bibitem{SS78}  Saulis, L. and Statulevi\v{c}ius, V. A., 1978. \textit{Limite theorems for large deviations}. Kluwer Academic Publishers.

}
\end{thebibliography}
\end{document}